\documentclass[12pt,reqno,a4paper]{amsart}
\usepackage{amsmath,amssymb,amsfonts,amsthm,amscd}
\usepackage[mathscr]{eucal}
\usepackage{hyperref}

\author{Ted Voronov}
\address{School of Mathematics, University of Manchester,
Oxford Road, Manchester, M13 9PL, United Kingdom}
\email{theodore.voronov@manchester.ac.uk}

\title[$Q$-manifolds and Mackenzie theory]{$Q$-manifolds and Mackenzie theory: an overview}

\thanks{\sc Based on the talk at the Erwin Schr\"{o}dinger Institute on 22 August (4 September) 2007}



\newtheorem*{them}{Theorem}
\newtheorem*{genprinc}{General principle}

\newtheorem{prop}{Proposition}[section]

\theoremstyle{definition}
\newtheorem{ex}{Example}[section]
\newtheorem{rem}{Remark}[section]

\def\co{\colon\thinspace}

\renewcommand{\div}{\mathop{\mathrm{div}}}

\DeclareMathOperator{\Vol}{Vol} 
 
 \DeclareMathOperator{\Vect}{\mathfrak{X}}

\newcommand{\der}[2]{{\frac{\partial {#1}}{\partial {#2}}}}

\newcommand{\fun}{C^{\infty}}

\newcommand{\w}{{\boldsymbol{w}}}

\renewcommand{\a}{\alpha}

\def\f{{\varphi}}

\newcommand{\h}{\eta}

\newcommand{\x}{{\xi}}

\newcommand{\ft}{{\tilde f}}

\newcommand{\jtt}{{\tilde \jmath}}

\newcommand{\ut}{{\tilde u}}

\begin{document}

\begin{abstract}
``Mackenzie theory" stands for the rich circle of notions that have been put forward by Kirill Mackenzie (solo or in collaboration): double structures such as double Lie groupoids and double Lie algebroids, Lie bialgebroids and their doubles, non-trivial dualities for double and multiple vector bundles, etc. ``$Q$-manifolds" are (super)manifolds --- I normally omit the prefix ``super" --- with a homological vector field, i.e., a self-commuting odd vector field. They may have an extra $\mathbb{Z}$-grading (called weight) not necessarily linked with the $\mathbb{Z}_2$-grading (parity). I will speak about double Lie algebroids (discovered by Mackenzie) and explain how this quite complicated fundamental notion is equivalent to a very simple one if the language of $Q$-manifolds is used. In particular, it shows how the two seemingly different notions of a ``Drinfeld double" of a Lie bialgebroid due to Mackenzie and Roytenberg respectively, turn out to be the same thing if properly understood.
\end{abstract}

\maketitle \tableofcontents

\section{Introduction}

This text is meant to be a brief overview of the topics announced in the title. It  does not contain new results (except probably for the remark concerning  $Q$-manifold homology, which we wish to elaborate elsewhere). The original exposition of Mackenzie's constructions is in his papers~\cite{mackenzie:secondorder1, mackenzie:bialg, mackenzie:secondorder2, mackenzie:doublealg, mackenzie:doublealg2, mackenzie:drinfeld,  mackenzie:sympldouble, mackenzie:notions, mackenzie:diffeomorphisms, mackenzie:duality}. See also the book~\cite{mackenzie:book2005}. There is a paper~\cite{mackenzie:ehresman}, which is now in the process of publication and which is a substantially reworked version of the earlier paper~\cite{mackenzie:notions}.
The main point of this text is to give an introduction to the approach to double Lie algebroids based on (graded) $Q$-manifolds. Details can be found in \cite{tv:double}, where the result (the equivalence between Mackenzie's double Lie algebroids and a certain class of $Q$-manifolds, see below) was formulated and proved for the first time. Some background information can also be found in~\cite{tv:graded}.

Double Lie algebroids first appeared as the infinitesimals  corresponding to the double Lie groupoids introduced by Mackenzie~\cite{mackenzie:secondorder1, mackenzie:secondorder2} as double objects in the sense  of Ehresmann (Lie groupoid objects in the category of Lie groupoids). The abstract notion came about later~\cite{mackenzie:doublealg, mackenzie:doublealg2, mackenzie:drinfeld} and turned out to be quite complicated in formulation. (One of the  reasons is that  a Lie algebroid is not defined diagrammatically, so it is not possible to follow a categorical approach, which works well for the groupoid case.) Although there was absolutely no doubt that this was the `right' notion --- one justification was in the appearance of double Lie algebroids as `Drinfeld doubles' of Lie bialgebroids~\cite{mackenzie:drinfeld} --- their application was somewhat hindered by the complexity of the definition. It had been of a considerable interest for experts to give an alternative simpler description for them. This was achieved in~\cite{tv:double}.
%There is also a work in progress~\cite{mackenzie:bidouble}.

This paper is an outcome of my participation in the program  on
``Poisson sigma models, Lie algebroids, deformations, and higher
analogues'' at the Erwin Schr\"odinger Institute in Vienna in
August--Sep\-tem\-ber 2007, where it was written. It is a pleasant
duty to express my gratitude to the Institute for the  hospitality
and the wonderful atmosphere for research and communications,  and
to the organizers of the program (Thomas Strobl,  Henrique Bursztyn
and Harald Grosse) for the invitation. The first announcement of the
results of~\cite{tv:double} was made in July 2006 at the XXV
Bia{\l}owie\.{z}a meeting.

\section{Lie algebroids and Lie bialgebroids}

First, let us recall some well-known facts concerning \textbf{Lie
algebroids}.

A \textit{Lie algebroid} is a vector bundle $E\to M$ over a manifold
$M$ endowed with a vector bundle map $a{\co} E\to TM$ (called the
\textit{anchor}) and a Lie algebra structure on the space of
sections $\fun(M,E)$ satisfying
\begin{align}%\label{algebroid}
  [u,fv]&={a(u)}f\,v+(-1)^{\ut\ft}f[u,v] \label{lr1}\\
  a([u,v])&= [a(u),a(v)]\label{lr2}
\end{align}
for all $u,v\in \fun(M,E)$ and all $f\in \fun(M)$.

This notion allows three other equivalent descriptions. Recall that a ${Q}$-manifold is a supermanifold endowed with a homological vector field (often denoted by $Q$). We shall also consider Poisson manifolds and odd Poisson (or Schouten) manifolds, sometimes using the abbreviations $P$- and $S$-manifolds respectively. ($S$ for Schouten.)

By a \textit{graded manifold} we shall understand a supermanifold with a privileged class of atlases where the coordinates are assigned weights taking values in $\mathbb{Z}$ and the coordinate transformations are polynomial in coordinates with nonzero weights   respecting the total weight. It is also assumed that the coordinates with nonzero weights run over the whole $\mathbb{R}$ (no restriction on range). Note that in general there is no relation between weight and parity. A natural example of a graded manifold is the total space of a vector bundle, for which we assume by default the following graded structure: the coordinates on the base have zero weight, the linear coordinates on fibers are assigned weight $1$. (Any graded manifold having only non-negative weights decomposes into a tower of affine fibrations, the first level being a vector bundle.) See more in~\cite{tv:graded}.

\begin{prop} There is a one-to-one correspondence between the following objects:
\begin{itemize}
  \item Lie algebroids;
  \item vector bundles with a Poisson bracket  of weight $-1$ on the total space;
  \item vector bundles with a Schouten bracket  of weight $-1$ on the total space;
  \item vector bundles with a homological vector field  of weight $+1$ on the total space.
\end{itemize}
\end{prop}

Indeed, consider a vector bundle $E\to M$ and its three neighbors: the dual bundle $E^*$, the opposite bundle $\Pi E$ and the antidual $\Pi E^*$. (Here $\Pi$ denotes the parity reversion functor.)

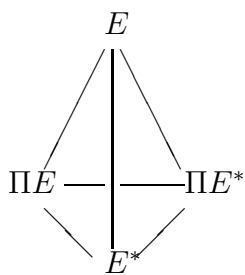
\begin{center}
{ \unitlength=3pt
\begin{picture}(0,35)
\put(0,30)
    {\begin{picture}(0,0)
    \put(0,0){$E$}
    \put(-12,-20){$\Pi E$}
    \put(10,-20){$\Pi E^*$}
    \put(0,-30){$E^*$}
    {\thicklines
    \put(1,-2){\line(0,-1){25}} }
    \put(-5,-19){\line(1,0){5}}
    \put(2,-19){\line(1,0){8}}
    \put(0,-2){\line(-1,-2){7.5}}
    \put(2,-2){\line(1,-2){7.5}}
    \put(-7.5,-22){\line(1,-1){6}}
    \put(10,-22){\line(-1,-1){6}}
    \end{picture}}
\end{picture}
}
\end{center}

A Lie algebroid structure for the bundle $E\to M$ is equivalent to either of the following: an even Poisson bracket  of weight $-1$ on the manifold $E^*$; an odd Poisson bracket  of weight $-1$ on the manifold $\Pi E^*$; and  a homological vector field  of weight $1$ on $\Pi E$. It is worth mentioning that a Lie algebroid structure on a bundle is introduced in terms of its sections, while each of the three other descriptions gives a structure on the total space.

Let me show the explicit formulas. Let $x^a$ be local coordinates on the base $M$ and $e_i$ make a local frame for $E\to M$. Then the corresponding coordinates on the fibers of $E^*$, $\Pi E$ and $\Pi E^*$ will be denoted by $u_i$, $\x^i$ and $\h_i$, respectively. (Here $\x^i e_i$ and similar expressions --- in this order --- are assumed to be invariant.) The correspondence between structures is as follows.

The anchor and Lie bracket for the basis sections of $E$:
\begin{equation*}
    a(e_i)=Q_i^a(x)\,\der{}{x^a}\,,
\end{equation*}
and
\begin{equation*}
    [e_i,e_j]=(-1)^{\jtt} Q_{ij}^k(x)\,e_k\,.
\end{equation*}

The (nonzero) even Poisson brackets  of the coordinates  on the manifold $E^*$ (I am skipping the precise signs):
\begin{equation*}
    \{x^a,u_i\}=\pm Q_i^a, \quad \{u_i,u_j\}=\pm Q_{ij}^k u_k\,.
\end{equation*}

Similarly, the odd brackets on $\Pi E^*$:
\begin{equation*}
    \{x^a,\h_i\}=\pm Q_i^a, \quad \{\h_i,\h_j\}=\pm Q_{ij}^k \h_k\,.
\end{equation*}

The homological field on $\Pi E$:
\begin{equation*}
    Q=\x^iQ_i^a(x)\,\der{}{x^a}+\frac{1}{2}\,\x^i\x^j Q_{ji}^k
    (x)\,\der{}{\x^k}\,.
\end{equation*}

(This is analogous to the similar statement for Lie algebras. For a Lie algebra $\mathfrak{g}$ we have the coalgebra $\mathfrak{g}^*$ with a linear Poisson structure,  the anticoalgebra $\Pi\mathfrak{g}^*$ with a linear Schouten structure, and the antialgebra $\Pi \mathfrak{g}$ with a quadratic homological field. Note that when we pass to vector bundles, we have to use the invariant notion of weight instead of speaking of something being `linear' or `quadratic' in coordinates.)

We  call a vector bundle with a homological vector field  of weight $+1$ on the total space, a \textit{Lie antialgebroid}.

(An arbitrary graded $Q$-manifold with $\w(Q)=1$, may therefore be called a \textit{generalized Lie antialgebroid}.)

\begin{rem}
If we take a Lie antialgebroid $\Pi E$ as a primary object, then the anchor and the Lie bracket for $E$ can be expressed by coordinate-free formulas using the derived bracket construction:
\begin{equation}\label{eq.auf}
    a(u)f:=\bigl[[Q,i(u)],f\bigr]
\end{equation}
and
\begin{equation} \label{eq.uv}
    i([u,v]):=(-1)^{\ut}\bigl[[Q,i(u)],i(v)\bigr].
\end{equation}
Here the map $i\co \fun(M, E)\to \Vect (\Pi E)$ has the following
appearance in coordinates:
\begin{equation}\label{eq.iu}
    i(u)=(-1)^{\ut}u^i(x)\der{}{\x^i}\,.
\end{equation}
(See, for example,~\cite{tv:graded}. See~\cite{yvette:derived,
yvette:derived2} for the derived bracket construction.
See~\cite{tv:higherder} for higher derived brackets.)
\end{rem}

Although the four descriptions are equivalent in the sense that they
all carry the same information, it turns out that the description in
the language of Lie antialgebroids is the most efficient of all. As
an example let is discuss the notion of a \textbf{Lie algebroid
morphism}. If  $E_1\to M_1$ and $E_2\to M_2$ are Lie algebroids over
different bases, the definition of a morphism is not obvious. One
first has to define a \textit{pull-back Lie algebroid} $\f^{!!}E_2$
w.r.t. a map of bases $\f\co M_1\to M_2$ as the `vector bundle'
$TM_1\times_{TM_2} E_2\to M_1$ with a certain anchor and Lie
bracket, which is not, strictly speaking, a genuine Lie algebroid,
because it is not in general a true vector bundle. Still, by this
the problem reduces to defining a morphism of Lie algebroids over
the same base $M_1$, which is straightforward. (The precise
definitions of the pull-back construction and of  Lie algebroid
morphisms over different bases can be found in the original
papers~\cite{mackenzie:noteonalg87, higgins-mackenzie:algebraic90}
and the book~\cite{mackenzie:book2005}.) The power of the language
of $Q$-manifolds can be illustrated by the following
\begin{prop} Consider Lie algebroids $E_1\to M_1$ and $E_2\to M_2$. Let $\Pi E_1$ and $\Pi E_2$ be the corresponding antialgebroids defined by the vector  fields $Q_1\in \Vect(\Pi E_1)$ and $Q_2\in \Vect(\Pi E_2)$ respectively.
A vector bundle map given by the horizontal
arrows of
\begin{equation*}
    \begin{CD} E_1 @>{\Phi}>> E_2 \\
                @VVV    @VVV\\
                M_1 @>{\varphi}>> M_2
    \end{CD}
\end{equation*}
is a morphism of Lie algebroids if and only if  the vector  fields $Q_1$ and $Q_2$ are
$\Phi^{\Pi}$-related, where $$\Phi^{\Pi}\co \Pi E_1\to \Pi E_2$$ is
the  induced map of the opposite vector bundles.
\end{prop}

(This characterization was given by Vaintrob~\cite{vaintrob:algebroids}.)

Shortly: a map of vector bundles is a Lie algebroid morphism if the induced map of antialgebroids is a morphism of $Q$-manifolds (which is much easier to handle).

Another illustration of the usefulness of the approach  based on
$Q$-manifolds  is an analog of \textbf{homology  for Lie
algebroids}. More generally, for a $Q$-manifold $M$, the
\textit{standard cochain complex} is $(C^{\infty}(M), Q)$. The
\textit{dual complex} or the \textit{chain complex} can be defined
as $(\Vol(M), L_Q)$.  Here $\Vol (M)$ stands for the Berezin volume
forms and $L_Q$, for the Lie derivative w.r.t. the vector field $Q$.
Note that `chains' are defined as cochains with certain
``dualizing'' coefficients. We are not looking at the
$\mathbb{Z}$-grading here; it can be taken care of when necessary.
There is a natural isomorphism (a ``Poincar{\'{e}} duality'')
between the cochain complex $(C^{\infty}(M), Q)$ and the chain
complex $(\Vol(M), L_Q)$ if there is an invariant non-vanishing
volume form $\rho$. One can easily see that the necessary and
sufficient condition for this is that the cohomology class
$[\div_{\boldsymbol{\rho}}Q]$ (which is independent of $\rho$) in
the standard complex equals zero. We call it the \textit{modular
class} of a $Q$-manifold. In the Lie algebroid case it is the
modular class introduced in~\cite{weinstein:evenslu99}. The point of
considering $(\Vol(M), L_Q)$ as the chain complex (i.e., giving the
homology of a $Q$-manifold) is that it has a natural pairing with
the cochain complex (provided certain conditions of compactness and
orientability are fulfilled, otherwise one has to consider
compactly-supported densities). It is not to be confused with the
mentioned  Poincar{\'{e}} duality. Most important, it possesses the
correct functorial behavior: for a map of $Q$-manifolds $\phi$ there
is a push-forward chain map $\phi_*$, hence there is an induced
push-forward map on homology. This is applicable to the particular
case of Lie algebroids, for which we obtain a homology theory with a
proper behavior under morphisms. For a Lie algebroid $E\to M$, the
chain complex $(\Vol(\Pi E), L_Q)$ coincides with the complex
$C(E,Q_E)$ introduced in a different language by Evens, Lu and
Weinstein~\cite{weinstein:evenslu99}. ($Q_E$ is their notation has
nothing to do with homological vector fields and stands for a
representation of $E$ defined in~\cite{weinstein:evenslu99}.) (For a
very particular case of the morphism $a\co E\to TM$   given by the
anchor, it is possible to show that the induced map of chain
complexes $a_{*}$ coincides with a map to differential forms on $M$
introduced in the top degree in~\cite{weinstein:evenslu99} in
connection with a version of  Poincar{\'{e}} duality. The fact that
it should be viewed as a special case of a push-forward, which is a
chain map in all degrees, was not known in the literature. We hope
to elaborate it elsewhere\,\footnote{This remark resulted from
discussions with Vladimir Roubtsov at ESI in August 2007.}.)

Now let us turn to  \textbf{Lie bialgebroids}. A   \textit{Lie bialgebroid}  over a supermanifold $M$ is a  Lie algebroid $E\to M$ such that the dual bundle $E^*\to M$ is also endowed with a structure of a Lie algebroid and a compatibility condition is satisfied.  There are several equivalent ways of expressing this compatibility condition. The original condition was given by Mackenzie and Xu~\cite{mackenzie:bialg}, who first introduced Lie bialgebroids; it was then replaced by a more efficient description due to Kosmann-Schwarzbach~\cite{yvette:exact}, which in its turn can be reformulated using $Q$-manifolds.
Consider  $\Pi E$. It  is a Lie antialgebroid because $E$ is a Lie algebroid, and it also carries a Schouten bracket of weight $-1$ because $E^*$ is a Lie algebroid. The condition is that the homological vector field is a derivation of the bracket. Shortly: $\Pi E$ is a $QS$-manifold of weights $(1,-1)$. One can show that this condition is symmetric in $E,E^*$ (see later), though it is not manifestly symmetric. Hence, Lie bialgebroids is a self-dual notion. They are often denoted as pairs $(E,E^*)$.

\section{Double and multiple Lie algebroids}

What is a double Lie algebroid? Such objects  naturally appear as
infinitesimals, i.e., by application of the Lie functor, for double
Lie groupoids (taken in the sense of Ehresmann, as groupoid objects
in the category of groupoids). Double Lie algebroids of double Lie
groupoids were introduced by Mackenzie
in~\cite{mackenzie:secondorder1}. The corresponding abstract
notion~\cite{mackenzie:secondorder2} appeared later and turned out,
in its original form, to be quite complicated and non-obvious. The
main reason for this, is that properties of brackets for Lie
algebroids are not expressed diagrammatically, so one cannot
approach double objects for them by methods of category theory as
one does for double groupoids. There is no doubt, however, that it
is a correct notion; this was shown, in particular, by the fact that
double Lie algebroids in the sense of Mackenzie naturally appear
also in the constructions of a `Drinfeld double' for Lie
bialgebroids~\cite{mackenzie:doublealg, mackenzie:drinfeld} (see
later). Only quite recently an alternative simplifying description
in terms of $Q$-manifolds was obtained~\cite{tv:double}, which also
helped to simplify, in the hindsight, Mackenzie's original approach
(see~\cite{mackenzie:ehresman}).

A double Lie algebroid    is, first of all,  a \textbf{double vector bundle}. Let us recall this notion. A \textit{double vector bundle} with base $M$ is a fiber bundle  $D\to M$ with a special structure. A local model for it (a trivial double vector bundle) has the appearance $U\times V_1\times V_2\times V_{12}$ where $U\subset M$ is an open set and $V_i, V_{ij}$ are vector spaces. Admissible transformations are maps $V_1\times V_2\times V_{12}\to V_1\times V_2\times V_{12}$, depending on points of $U$ as  parameters,  which for each $V_i$ are linear  and for $V_{12}$ linear in $V_{12}$ plus an extra  term bilinear in $V_1\times V_2$. In other words, if we use linear coordinates $u^i$ for $V_1$, $w^{\a}$ for $V_2$, and $z^{\mu}$ for $V_{12}$, we have a transformation law of the form
\begin{align}
    u^i&=u^{i'}T_{i'}{}^{i}(x'),  \label{eq.lawforu}\\
    w^{\a}&=w^{\a'}T_{\a'}{}^{\a}(x'), \label{eq.lawforw} \\
    z^{\mu}&=z^{\mu'}T_{\mu'}{}^{\mu}(x')+w^{\a'}u^{i'}T_{i'\a'}{}^{\mu}(x')\,.\label{eq.lawforz}
\end{align}
(all coefficients are functions of $x\in U$). In particular, it follows that there is a diagram
\begin{equation} \label{eq.dvb2}
    \begin{CD} D@>>> B\\
                @VVV  @VVV \\
                A@>>>M
    \end{CD}
\end{equation}
for which each side is a vector bundle. Here $V_1$ is the standard fiber for $A\to M$; $V_2$, for $B\to M$; $V_1\times V_{12}$, for $D\to B$; and $V_2\times V_{12}$, for $D\to A$. There is also a vector bundle $K\to M$ with the standard fiber $V_{12}$. It is called the \textit{core} of the double vector bundle $D\to M$.

(There is an immediate generalization giving the notion of a \textit{$k$-fold vector bundle}. Such a bundle has the standard fiber of the form $\prod V_i\times \prod_{i<j}V_{ij}\times \ldots \times V_{12\ldots k}$ and the transformation law similar to the above. By picking a subset $i_1<\ldots<i_l$ one obtains a \textit{face}, which is an $l$-fold vector bundle. A $1$-fold vector bundle is an ordinary vector bundle.)

From the definition it is clear that the total space of a double
vector bundle $D\to M$  is naturally a bigraded manifold, and the
total space of a  $k$-fold vector bundle $P\to M$ is a $k$-graded
manifold. One can speak about the total weight as the sum of all
(partial) weights.

There are constructions that naturally lead to multiple vector bundles.

\begin{ex} If $E\to M$ is a vector bundle, then the tangent $TE$ has the structure of a double vector bundle:
\begin{equation*} %\label{eq.dvb2}
    \begin{CD} TE@>>> E\\
                @V{Tp}VV  @VV{p}V \\
                TM@>>>M
    \end{CD}
\end{equation*}
To check that it is indeed a double vector bundle, one has to write down the changes of coordinates. If
\begin{equation*}
    v^i=v^{i'}T_{i'}^{i}
\end{equation*}
is the transformation law for the fiber coordinates on $E$, then for $TE$ we obtain the fiber coordinates $v^i,\dot x^a, \dot v^i$ where $x^a$ are local coordinates on the base $M$, and clearly
\begin{align}
    \dot x^a&=\dot x^{a'}\der{x^a}{x^{a'}}\\
    \dot v^i&= \dot v^{i'}T_{i'}^{i}+ v^{i'} \dot x^{a'}\der{T_{i'}^{i}}{x^{a'}}\,,\label{eq.dotvi}
\end{align}
which satisfies the definition of a double vector bundle. Here as the core we have a vector bundle with the transition functions obtained from~\eqref{eq.dotvi} by setting $v^i$ and $\dot x^a$ to zero; therefore it is a copy of $E\to M$.
\end{ex}

\begin{ex} The cotangent  $T^*E$, for a vector bundle $E\to M$, is also a double vector bundle. Here we have the diagram
\begin{equation*} %\label{eq.dvb2}
    \begin{CD} T^*E@>>> E\\
                @VVV  @VVV \\
                E^*@>>>M
    \end{CD}
\end{equation*}
The core bundle in this case is $T^*M\to M$. We shall return to this example later.
\end{ex}

(See book~\cite{mackenzie:book2005}.)

Both examples can be generalized to the case when we already start from a multiple vector bundle; taking the tangent or cotangent gives a  multiple vector bundle of the multiplicity increased by one.

Now the definition of a \textit{double Lie algebroid} with base $M$ is as follows: it is a double vector bundle
\begin{equation} \label{eq.dvb3}
    \begin{CD} D@>>> B\\
                @VVV  @VVV \\
                A@>>>M
    \end{CD}
\end{equation}
such that each side (which is a vector bundle) carries the structure of a Lie algebroid and certain compatibility conditions are satisfied. The main problem happens to be to formulate these compatibility conditions. To this we are going to proceed, but we will have to make a further digression to constructions with multiple vector bundles.

Double (and multiple) vector bundles possess a non-trivial \textbf{duality theory}. Let $D\to M$ be a double vector bundle. There is a diagram of ordinary vector bundles~\eqref{eq.dvb2}. Consider $D$  as a vector bundle over $A$ and take the dual, which we denote $D^*_A$ explicitly indicating the base.

\begin{prop} $D^*_A$ is again a double vector bundle over $M$, given by a square
\begin{equation*} %\label{eq.dvb2}
    \begin{CD} D^*_A@>>> K^*\\
                @VVV  @VVV \\
                A@>>>M
    \end{CD}
\end{equation*}
where $K^*\to M$ is the dual to the core $K\to M$.
\end{prop}
(Of course, the same is true if we interchange $A$ and $B$.)

This can be explained using coordinates (or local trivializations) as follows.  In the notation $x^a,u^i,w^{\a},z^{\mu}$ for coordinates on $D$ used above, for $D^*_A$ we must have local coordinates $x^a,u^i,w_{\a},z_{\mu}$ so that the form $w^{\a}w_{\a}+z^{\mu}z_{\mu}$ be invariant. The invariance condition together with the transformation law~\eqref{eq.lawforw} and~\eqref{eq.lawforz}  imply the transformation law
\begin{align}
    w_{\a'}&=T_{\a'}{}^{\a}w_{\a}+ u^{i'}T_{i'\a'}^{\mu}z_{\mu}, \label{eq.lawforwniz} \\
    z_{\mu'}&=T_{\mu'}{}^{\mu}z_{\mu}\,.\label{eq.lawforzniz}
\end{align}
We see that a double vector bundle is indeed obtained, and the new core is $B^*\to M$, the dual  for $B\to M$. If we apply duality to the double vector bundle $D^*_A$, then, depending in which direction we dualize, we either come back to $D$ or arrive at the other dual $D^*_B$. Hence \textit{there are three double vector bundles in duality}, making a `corner':
\begin{equation}
    \begin{picture}(300,80)(0,65)
    \put(90,90){${\begin{CD} {D^*_A}  @>>>   K^*\\
                @VVV  @VVV \\
               A@>>>M\end{CD}}$}
    \put(120,110){$\begin{CD} {} @.   D^*_B\\
                @.  @VVV \\
               D@>>>B\end{CD}$}
    {\put(120,84){\vector(-2,-1){17}}}
    {\put(180,84){\vector(-2,-1){17}}}
    %\put(132,130){\vector(-2,-1){18}}
    {\put(180,130){\vector(-2,-1){17}}}
    \end{picture}
\end{equation}
This amazing duality between the `horizontal' and `vertical'  duals
of a double vector bundles was discovered by
Mackenzie~\cite{mackenzie:sympldouble} and independently by
Konieczna--Urba\'{n}ski~\cite{konieczna:duality}. The canonical
pairing between $D^*_A$ and $D^*_B$ as vector bundles over $K^*$ is
given by the invariant bilinear form
\begin{equation}\label{eq:pairing}
    u^iu_i-w^{\a}w_{\a}\,.
\end{equation}
(Recall that $u^i,w_{\a}$ and $u_i,w^{\a}$ are precisely the fiber
coordinates  on $D^*_A$ and $D^*_B$ over $K^*$, so that
$x^a,z_{\mu}$ are the base coordinates.) The main statement is the
invariance of the form~\eqref{eq:pairing}, for which the minus sign
is absolutely essential (the two terms being not separately
invariant).

Similar statements hold  for triple and multiple vector bundles. (See Mackenzie~\cite{mackenzie:duality} for a detailed analysis.)

Besides duality functors another type of operation that can be applied to multiple vector bundles is changing parity in fibers. As the parity reversion  can be applied only to linear objects, there are $k$ partial \textbf{parity reversion} functors $\Pi_i$ for a $k$-fold vector bundle. One can show that
\begin{equation*}
    \Pi_i\Pi_j=\Pi_j\Pi_i\,,
\end{equation*}
meaning natural isomorphism.

Now we can proceed to double algebroids. Consider the double vector bundle given by~\eqref{eq.dvb3}. The obvious part of the definition is that all sides should carry a Lie algebroid structure. Now, the original definition due to Mackenzie contained three extra conditions: the so-called \textbf{Conditions I}, \textbf{II}, and \textbf{III}, which are rather complicated in their formulation. We shall not reproduce them here, because it would require a different level of technical detail than we are using. Instead we shall give an informal idea of their content.

Condition I is the easiest for understanding. It basically requires
that each of the algebroid structures respects the linear structure
in the other direction. It turns out, in particular, that the Lie
algebroid structures on $D\to A$ and $D\to B$ induce those on $B$
and $A$. Besides this, if we consider homological vector fields
$Q_{DA}$ on $\Pi_AD$ and $Q_{DB}$ on $\Pi_BD$ corresponding to these
Lie algebroid structures, they must have zero weights in the other
direction, i.e.,
\begin{equation*}
    \w(Q_{DA})=(1,0), \quad \w(Q_{DB})=(0,1)\,.
\end{equation*}
Note that weight zero means that a vector field generates linear transformations. This, in particular, allows to dualize and to reverse parity.

Now instead of giving a technical description of Mackenzie's conditions II and III, let us turn to what I call a \textbf{big picture}. The lesson that we learned from considering ordinary Lie algebroids is that this notion has different (and equivalent) manifestations, which we see by looking at all the neighbors of a given vector bundle. Let us do the same for the double vector bundle~\eqref{eq.dvb3}. In the list below we indicate which structures correspond to the two Lie algebroid structures on $D\to A$ and $D\to B$. Altogether there are $12$ different neighbors (counting the original double vector bundle).

\begin{equation} \label{eq.lieandpoiss}
    \begin{CD} D^{*}_{A}@>>>  K^*\\
                @VVV  @VVV \\
                A@>>>M
    \end{CD} \qquad \text{\small P. bracket on $D^*_A$, $(-1,-1)$; L. algd. on $D^*_A\to K^*$}
\end{equation}

\begin{equation} \label{eq.poissandlie}
    \begin{CD} D^{*}_{B}@>>>  B\\
                @VVV  @VVV \\
                K^*@>>>M
    \end{CD} \qquad \text{\small P. bracket on $D^*_B$, $(-1,-1)$; L. algd. on $D^*_B\to K^*$}
\end{equation}

\begin{equation} \label{eq.dv2and3}
    \begin{CD} \Pi_AD@>>> \Pi B\\
                @VVV  @VVV \\
                A@>>>M
    \end{CD} \qquad \text{\small H. v. field on $\Pi_AD$, $(1,0)$; L. algd. on $\Pi_AD\to \Pi B$}
\end{equation}

\begin{equation}
    \begin{CD} \Pi_BD@>>>  B\\
                @VVV  @VVV \\
                \Pi A@>>>M
    \end{CD} \qquad \text{\small H. v. field on $\Pi_BD$, $(0,1)$; L. algd. on $\Pi_BD\to \Pi A$}
\end{equation}

\begin{equation} \label{eq.pikvadratt}%\label{eq.dv4}
    \begin{CD} \Pi^2D@>>>  \Pi B\\
                @VVV  @VVV \\
                \Pi A@>>>M
    \end{CD} \qquad \text{\small Two h. v. fields on $\Pi^2D$, weights $(0,1)$ and $(1,0)$}
\end{equation}

\begin{equation} \label{eq.dv7}
    \begin{CD} \Pi_AD^{*}_{A}@>>>  \Pi K^*\\
                @VVV  @VVV \\
                A@>>>M
    \end{CD} \qquad \text{\small S. bracket on $\Pi_AD^{*}_{A}$, $(-1,-1)$; L. algd. on $\Pi_AD^{*}_{A}\to \Pi K^*$}
\end{equation}

\begin{equation} \label{eq.dv11}
    \begin{CD} \Pi_{B}D^{*}_{B}@>>>   B\\
                @VVV  @VVV \\
               \Pi K^*@>>>M
    \end{CD} \qquad \text{\small S. bracket on $\Pi_BD^{*}_{B}$, $(-1,-1)$; L. algd. on $\Pi_BD^{*}_{B}\to \Pi K^*$}
\end{equation}

\begin{equation} \label{eq.piduala}%\label{eq.dv8}
    \begin{CD} \Pi_{K^*}D^{*}_{A}@>>>   K^*\\
                @VVV  @VVV \\
                \Pi A@>>>M
    \end{CD} \qquad \text{\small S. bracket, $(-1,-1)$, and h. v. field, $(1,0)$, on $\Pi_{K^*}D^{*}_{A}$}
\end{equation}

\begin{equation} \label{eq.pidualb}%\label{eq.dv10}
    \begin{CD} \Pi_{K^*}D^{*}_{B}@>>>  \Pi B\\
                @VVV  @VVV \\
                K^*@>>>M
    \end{CD} \qquad \text{\small S. bracket, $(-1,-1)$, and h. v. field, $(0,1)$, on $\Pi_{K^*}D^{*}_{B}$}
\end{equation}

\begin{equation} \label{eq.pikvadratduala}%\label{eq.dv9}
    \begin{CD} \Pi^2 D^{*}_{A}@>>>   \Pi K^*\\
                @VVV  @VVV \\
                \Pi A@>>>M
    \end{CD} \qquad \text{\small P. bracket, $(-1,-1)$, and h. v. field, $(1,0)$, on $\Pi^2 D^{*}_{A}$}
\end{equation}

\begin{equation} \label{eq.pikvadratdualb}%\label{eq.dv12}
    \begin{CD} \Pi^2 D^{*}_{B}@>>>   \Pi B\\
                @VVV  @VVV \\
                \Pi K^*@>>>M
    \end{CD} \qquad \text{\small P. bracket, $(-1,-1)$, and h. v. field, $(0,1)$, on $\Pi^2 D^{*}_{B}$}
\end{equation}

Now, the general philosophy is: if we do not know how to define
compatibility of the two structures for some of these pictures, look
at the neighbors for which a compatibility condition comes about
naturally. In our list there are the exactly five double vector
bundles given by diagrams~\eqref{eq.pikvadratt}, \eqref{eq.piduala},
\eqref{eq.pidualb},  \eqref{eq.pikvadratduala}, and
\eqref{eq.pikvadratdualb}, where both structures are defined on the
total space. It is absolutely clear which condition should be
considered in each of these cases: the commutativity of the vector
fields for \eqref{eq.pikvadratt} and the derivation property w.r.t.
the bracket for \eqref{eq.piduala},  \eqref{eq.pidualb},
\eqref{eq.pikvadratduala}, and \eqref{eq.pikvadratdualb}.

\begin{prop} \label{prop.bialg}
The compatibility conditions for \eqref{eq.piduala},
\eqref{eq.pidualb}, \eqref{eq.pikvadratduala} and
\eqref{eq.pikvadratdualb} are equivalent, and are different ways of
saying that $(D^{*}_{A}, D^{*}_{B})$ is a Lie bialgebroid over
$K^*$.
\end{prop}
\begin{proof} Indeed, for any Lie bialgebroid $(E,E^*)$ the compatibility can
be stated in terms of either $E$ or $E^*$ (as a $QS$-structure on
either $\Pi E$ or $\Pi E^*$, respectively). This corresponds
to~\eqref{eq.piduala} or \eqref{eq.pidualb} in our case. In our
special case  there is also an extra option of changing parity in
the second direction, which adds~\eqref{eq.pikvadratduala} and
\eqref{eq.pikvadratdualb} to the picture.
\end{proof}

On the other hand, to say  that $(D^{*}_{A}, D^{*}_{B})$ is a Lie
bialgebroid over $K^*$ is a natural compatibility condition for
\eqref{eq.lieandpoiss} and \eqref{eq.poissandlie}. And it is
precisely Mackenzie's Condition III. It was found
in~\cite{tv:double} and then proved by Mackenzie in his own
framework~\cite{mackenzie:ehresman} that his original Condition II
is subsumed by Condition III (though it is not obvious).

\begin{them} A double vector bundle $D\to M$ is a double Lie algebroid  if and only if on the double vector bundle $\Pi^2D\to M$ the two induced homological vector fields (corresponding to the two Lie algebroid structures) commute.
\end{them}
\begin{proof} Consider one of the manifestations of the bialgebroid
condition, say, for concreteness,~\eqref{eq.piduala}. The derivation
property means that the flow of the vector field preserves the
bracket. On the other hand, the commutativity condition
for~\eqref{eq.pikvadratt} means that the flow of one field preserves
the other. Now the claim follows from  functoriality: notice that a
linear transformation preserves a Lie bracket if and only if the
adjoint map preserves the corresponding linear Poisson bracket and
if and only if the `$\Pi$-symmetric' map preserves the corresponding
homological vector field.
\end{proof}

We may define a \textit{double Lie antialgebroid} as  a double
vector bundle endowed with two commuting homological vector fields
of weights $(1,0)$ and $(0,1)$. In a similar way one defines a
\textit{$k$-fold Lie antialgebroid}. It is clear that the structures
of  double Lie algebroids and double Lie antialgebroids are
equivalent. This also gives an efficient general notion of a
\textit{multiple Lie algebroid}.

\section{Example: Drinfeld double for Lie bialgebroids}

Recall that Drinfeld's \textit{classical double} of a Lie bialgebra
is again a Lie bialgebra with ``good'' properties. An analog of this
construction for Lie algebroids turned out to be a puzzle.  Three
constructions of a `double' of a Lie bialgebroid have been
suggested. Suppose $(E,E^*)$ is a Lie bialgebroid over a base $M$.
Liu, Weinstein and Xu~\cite{weinstein:liuxu} suggested to consider
as its double a structure of a Courant algebroid on the direct sum
$E\oplus E^*$. Mackenzie in~\cite{mackenzie:doublealg,
mackenzie:doublealg2, mackenzie:drinfeld, mackenzie:notions} and
Roytenberg in~\cite{roytenberg:thesis} suggested two different
constructions based on cotangent bundles. Though they look very
different (in particular, Roytenberg's double is a supermanifold,
and Mackenzie stays in the classical world), we shall  show now that
they are essentially the same.

Roytenberg previously showed~\cite{roytenberg:thesis} that the
Liu--Weinstein--Xu double is recovered from his own construction as
a derived bracket, generalizing the results  of
C.~Roger~\cite{roger:1991} and
Y.~Kosmann-Schwarzbach~\cite{yvette:jacobian, yvette:derived} for
Lie bialgebras. Therefore, proving that the Mackenzie and Roytenberg
pictures  are equivalent or, actually,   the same, if understood
properly, shows conclusively that this `cotangent double' is
fundamental, and should be regarded as the correct extension of
Drinfeld's double of Lie bialgebras to Lie bialgebroids.

Both Roytenberg's and Mackenzie's construction use the statement
that the cotangent bundles of dual vector bundles are isomorphic
(\cite{mackenzie:bialg}, an extension of ~\cite{tulczyjew:1977}; see
also~\cite{mackenzie:diffeomorphisms}, \cite{roytenberg:thesis},
\cite{tv:graded}). Hence there is a double vector bundle
\begin{equation} \label{eq.mackdouble}
    \begin{CD} T^*E=T^*E^*@>>>   E^*\\
                @VVV  @VVV \\
               E@>>>M
    \end{CD}
\end{equation}
Mackenzie shows that it is a double Lie algebroid. He calls it the
\textit{cotangent double} of a Lie bialgebroid $(E,E^*)$. Note that
the canonical symplectic structure on $T^*E$ corresponds to the
invariant scalar product on Drinfeld's double $\mathfrak d
(\mathfrak b)=\mathfrak b\oplus \mathfrak b^*$ of a Lie bialgebra
$\mathfrak b$.

On the other hand, Roytenberg uses the description of Lie algebroids
via homological vector fields. He considers the double vector bundle
\begin{equation} \label{eq.roytdouble}
    \begin{CD} T^*\Pi E=T^*\Pi E^*@>>>   \Pi E^*\\
                @VVV  @VVV \\
               \Pi E@>>>M
    \end{CD}
\end{equation}
and   homological vector fields  $Q_E\in\Vect (\Pi E)$ and
$Q_{E^*}\in\Vect (\Pi E^*)$ defining Lie algebroid structures on
$E\to M$ and $E^*\to M$, respectively. Recall that  vector fields on
a manifold correspond to fiberwise linear functions (Hamiltonians)
on the cotangent bundle so that the commutator  maps to the Poisson
bracket. Denote the functions corresponding to $Q_E$ and $Q_{E^*}$
by $H_E$ and $H_{E^*}$, respectively. Roytenberg shows that under
the natural symplectomorphism $T^*\Pi E=T^*\Pi E^*$ the linear
function $H_{E^*}$ on $T^*\Pi E^*$ corresponding to the vector field
$Q_{E^*}$ transforms precisely into the fiberwise quadratic function
$S_E$ on $T^*\Pi E$  specifying the Schouten bracket on $\Pi E$
induced by the Lie structure on $E^*$. The derivation property of
$Q_E$ w.r.t. the Schouten bracket on $\Pi E$ is one of the
equivalent definitions of a Lie bialgebroid~\cite{yvette:exact}, and
the most convenient. Hence, Roytenberg's statement means that it is
also equivalent to the commutativity of the Hamiltonians $H_E$ and
$H_{E^*}$ under the canonical Poisson bracket. They generate
commuting homological vector fields $X_{H_E}$ and $X_{H_{E^*}}$ on
the cotangent bundle $T^*\Pi E$. In our language, $X_{H_E}$ and
$X_{H_{E^*}}$ make~\eqref{eq.roytdouble} a Lie antialgebroid. (One
can see that the conditions for weights are satisfied.)
Roytenberg~\cite{roytenberg:thesis} calls the supermanifold $T^*\Pi
E=T^*\Pi E^*$ together with the homological vector field
$Q=X_{H_E}+X_{H_{E^*}}$ on it, the \textit{Drinfeld double} of
$(E,E^*)$.

If we slightly refine Roytenberg's picture,  considering the double
Lie antialgebroid given by $X_{H_E}$ and $X_{H_{E^*}}$ rather than a
single $Q$-manifold, we can immediately see that by our
theorem  his picture becomes identical to that of
Mackenzie.

Indeed, apply the complete reversion of parity
to~\eqref{eq.mackdouble}. Notice that $\Pi^2 T^*E=\Pi^2 T^*E^*$
coincides with $T^* \Pi E=T^* \Pi E^*$ (easily checked in
coordinates). By the theorem, the double vector
bundle~\eqref{eq.mackdouble} is a double Lie algebroid if and only
if the corresponding double vector bundle
\begin{equation} \label{eq.mackdoublepikv}
    \begin{CD} \Pi^2T^*E=\Pi^2T^*E^*@>>>   \Pi E^*\\
                @VVV  @VVV \\
               \Pi E@>>>M
    \end{CD}
\end{equation}
which is identical with~\eqref{eq.roytdouble}, is a double Lie
antialgebroid. It remains to identify the respective homological
vector fields on the ultimate total space, which can be achieved by
a direct inspection.

We have arrived at the following statement.
\begin{prop} Roytenberg's and Mackenzie's pictures give the same notion of a double of a Lie
bialgebroid (up to a change of parity).
\end{prop}

We can now identify the two constructions and speak simply of the
(cotangent) \textit{double} of a Lie bialgebroid  as a double Lie
algebroid, most efficiently described in the ``anti-'' language of
diagram ~\eqref{eq.roytdouble}.

\textbf{More on doubles.}

\smallskip
Recall that Drinfeld's classical double of a Lie bialgebra is not
just a Lie algebra, but also a coalgebra, and in fact a Lie
bialgebra again. This gives a direction in which to look in the case
of Lie bialgebroids. Note that this second structure (for doubles of
Lie bialgebroids) has not been discovered previously.

However, for many people including the author there was absolutely
no doubt that such a structure exists. The following conjectured
statement was put down in my notes of 2002 as a guideline for (not
yet obtained at that time) alternative description of double Lie
algebroids.

\begin{genprinc}
Taking the double of an $n$-fold Lie bialgebroid should give  an
$(n+1)$-fold Lie bialgebroid, with an additional property, such as a
symplectic structure.
\end{genprinc}

Of course it involved new notions yet to be defined.   As for ``double Lie \textbf{bi}algebroids''
(or ``bi- double Lie algebroids''), and the further multiple ``bi-''
case, this  notion is  properly defined in our joint work with
Kirill Mackenzie, and is a subject of our forthcoming
paper~\cite{mackenzie:bidouble}, where precise definitions and
statements can be found.

%\bibliographystyle{plain}
%\bibliography{c:/Ted/Bibliography/geometry}
%\end{document}

\end{document}